\documentclass[12pt,final]{article}
\usepackage[english]{babel}
\usepackage{amssymb, amsmath, epsfig, graphicx, color, textcomp}

\def\Cset{\mathbb{C}}
\def\Nset{\mathbb{N}}

\def\Rset{\mathbb{R}}

\def\grad{\mbox{grad }}

\def\curl{\mbox{curl }}

\def\pro{\bf Proof: \rm}

\def\endbew{\hfill$\Box$\medskip}
\def\endpro{\endbew}

\def\clK2{\overline{K}_2}

\def\Id{$\Id$}

\renewcommand{\it}{k}
\newcommand{\It}{K}
\newcommand{\cgit}{l}
\newcommand{\dimX}{M}
\newcommand{\dimY}{N}
\newcommand{\calJ}{\mathcal{J}}
\newcommand{\El}{\mathbf{E}} 
\newcommand{\Eli}{\mathbf{E}^{\mathrm{i}}} 
\newcommand{\Els}{\mathbf{E}^{\mathrm{s}}} 
\newcommand{\Elinf}{\mathbf{E}^{\infty}} 
\newcommand{\bfd}{\mathbf{d}}
\newcommand{\bfp}{\mathbf{p}}
\newcommand{\bfr}{\mathbf{r}}
\newcommand{\bfx}{\mathbf{x}}
\newcommand{\bfy}{\mathbf{y}}
\newcommand{\supp}{\,\mathrm{supp}\,}

\newcommand{\D}{\,{\rm d}}

\newcommand{\bl}[1]{{\bf #1}}
\newcommand{\G}{\bl{G}}

\newcommand{\M}{\bl{M}}

\begin{document}

\newtheorem{Def}{Definition}
\newtheorem{tho}[Def]{Theorem}
\newtheorem{lem}[Def]{Lemma}
\newtheorem{cor}[Def]{Corollary}
\newtheorem{prop}[Def]{Proposition}
\newtheorem{example}[Def]{Example}
\newtheorem{assumption}[Def]{Assumption}
\newtheorem{rem}[Def]{Remark}
\newtheorem{notation}[Def]{Notation}
\newtheorem{alg}[Def]{Algorithm}

\title{Acceleration techniques for regularized Newton methods applied to 
electromagnetic inverse medium scattering problems}
\author{\Large{Thorsten Hohage%
\thanks{Institut f\"ur Numerische und Angewandte Mathematik, Universit\"at G\"ottingen
Lotzestr. 16--18, D-37083 G\"ottingen, Germany. E-mail: hohage@math.uni-goettingen.de}} 
and Stefan Langer\thanks{Deutsches Zentrum f\"ur Luft- und Raumfahrt e. V.
in der Helmholtz-Gemeinschaft, Institut f\"ur Aerodynamik und
Str\"omungstechnik, 
Lilienthalplatz 7,
38108 Braunschweig,
Germany. Tel.: +49 531 2953615; fax: +49 531 2952914. E-mail:
stefan.langer@dlr.de}}    
\date{\today}

\maketitle
 
%
%

\begin{abstract} 
We study the construction and updating of spectral preconditioners for´
regularized Newton methods and their application to electromagnetic
inverse medium scattering problems. Moreover, we show how a Lepski{\u\i}-type
stopping rule can be implemented efficiently for these methods.
In numerical examples, the proposed method compares favorably with other
iterative regularization method in terms of work-precision diagrams
for exact data. For data perturbed by random noise, the Lepski{\u\i}-type
stopping rule performs considerably better than the commonly used
discrepancy principle. 
\end{abstract}

\section[Introduction]{Introduction}
In this paper we study the efficient numerical solution of an inverse scattering 
problem for time harmonic electromagnetic waves. The forward problem  
is essentially described by the time-harmonic Maxwell equations 
\[
\curl\curl \El(\bfr) -\kappa^2 n(\bfr)^2\El(\bfr) = 0
\]
for the electric field $\El$. Our aim is to reconstruct a local inhomogeneity of 
the refractive index $n$ of a medium, given far field measurements for many incident waves.
A more detailed discussion of the forward problem is given in \S \ref{sec:forward}. 

After discretization the inverse problem is described by a nonlinear, ill-conditioned
system of equations $\bl{F}(\bl{x})=\bl{y}$ with a function $\bl{F}:D(\bl{F})\subset\Rset^{\dimX}\to \Rset^\dimY$, which is infinitely smooth on the subset $D(\bl{F})\subset\Rset^{\dimX}$ 
where it is defined. Since the system is highly ill-conditioned, we have consider 
the effects of data noise. Here we assume an additive noise model for the
observe data $\bl{y}^{\rm obs}$:
\begin{equation}\label{eq:noise_model}
\bl{y}^{\rm obs} =\bl{F}(\bl{x}) + \bl{\epsilon}
\end{equation}
The noise vector $\bl{\epsilon}$ is assumed to be a vector of random variables
with known finite covariance matrix and a known bound on the expectation
$\|\mathbb{E} \bl{\epsilon}\| \leq \delta$.


In this article we contribute to preconditioning techniques for the
\emph{Levenberg-Marquardt algorithm} and the \emph{iteratively regularized
Gauss-Newton method} (IRGNM). These methods are obtained by applying
Tikhonov regularization with some an initial guess $\bl{b}_{\it}$ and a 
regularization 
parameter $\gamma_{\it}$ to the Newton equations $\bl{A}_{\it}\bl{h}_{\it} 
= \bl{y}^{\rm obs}-\bl{F}(\bl{x}_{\it})$. Here $\bl{A}_{\it}:=\bl{F}'[\bl{x}_{\it}]\in \Rset^{\dimY\times \dimX}$
denotes the Jacobian of $\bl{F}$ at $\bl{x}_{\it}$. 
This leads to normal equations of the form 
\begin{equation}\label{260201}
\bl{G}_{\it}^{\top}\bl{G}_{\it}\bl{h}_{\it} = \bl{G}_{\it}^{\top}\bl{g}_{\it}
\end{equation}
with
\[
\bl{G}_{\it} := \left (\begin{array}{c} 
\bl{A}_{\it}\\ 
\sqrt{\gamma_{\it}}\bl{I}\end{array}\right ) \in\Rset^{(\dimY+\dimX)\times \dimX}\qquad\mbox{and}\qquad
\bl{g}_{\it} := \left (\begin{array}{c}
\bl{y}^{\rm obs} - \bl{F}(\bl{x}_{\it})\\ 
\sqrt{\gamma_{\it}}\bl{b}_{\it} \end{array}
\right). 
\]
The choice~$\bl{b}_{\it}=0$ corresponds to the Levenberg-Marquardt algorithm and
the choice $\bl{b}_{\it} = \bl{x}_0-\bl{x}_{\it}$ to  the IRGNM. 
As opposed to the Levenberg-Marquardt algorithm as used in optimization
we simply choose the regularization parameter~$\gamma_{\it}$ of the form
\begin{equation}\label{lch257}
\gamma_{\it} = \gamma_0\gamma^{-\it}\qquad \mbox{with } \gamma > 1.
\end{equation}
Convergence and convergence rates of the IRGNM in an infinite dimensional
setting have been studied first in \cite{bakush:92,BNS:97,hohage:97}. 
For further references and results including a convergence analysis of 
Levenberg-Marquardt algorithm we refer to the monographs
\cite{BaKo:04,KNS:08}.

As an alternative, Hanke \cite{hanke:97b} suggested to apply the conjugate 
gradient (CG) method the normal equation 
$\bl{A}_{\it}^{\top}\bl{A}_{\it}\bl{h}_{\it} 
= \bl{A}_{\it}^{\top}(\bl{y}^{\rm obs}-\bl{F}(\bl{x}_{\it}))$
and use the regularizing properties of the CG method applied
to the normal equation with early stopping. This is referred 
to as Newton-CG method. 
Regularized Newton methods with inner iterative regularization
methods have also been studied by Rieder \cite{rieder:99,rieder:05}.
Finally, applying a gradient method to the functional
$\bl{x}\mapsto \frac{\mu}{2}\|\bl{F}(\bl{x})-\bl{y}^{\rm obs}\|_2^2$ leads to the
nonlinear Landweber iteration $\bl{x}_{\it+1} :=\bl{x}_{\it}
-\mu \bl{A}_{\it}^{\top}(\bl{F}(\bl{x}_{\it})-\bl{y}^{\rm obs})$ first studied
in \cite{HNS:95}. For an overview on iterative regularization methods
for nonlinear ill-posed problems we refer to \cite{KNS:08}. 

A continuation method for inverse
electromagnetic medium scattering problems with multi-frequency
data has been studied in \cite{BL:05}.
For an overview on level set methods for inverse scattering problems
we refer to \cite{DL:06,DL:09}. 

\medskip
For the inverse electromagnetic scattering problem studied in this
paper the evaluation of $\bf{F}$ and one row of its Jacobian $\bf{A}_{\it}$
is very expensive and involves the solution of a three-dimensional
forward scattering problem for many incident waves. Therefore, 
a computation of the full Jacobian is not reasonable, and
regularization method for the inverse problem should only access
$\bf{A}_{\it}$ via matrix-vector multiplications 
$\bl{v}\mapsto\bl{A}_{\it}\bl{v}$ and 
$\bl{g}\mapsto \bl{A}_{\it}^{\top}\bl{g}$. 
Hence, from the methods discussed above only Landweber iteration and
Newton-CG can be implemented directly. However, the convergence of
Landweber iteration is known to be very slow, which is confirmed
by our numerical experiments reported in \S \ref{sec:num}.
Preconditioning techniques for Landweber iteration have been studied
in \cite{EN:05}, but it is not clear how to apply these techniques
to inverse electromagnetic medium scattering problems since the
operator does not act in Hilbert scales. To use the IRGNM and
Levenberg-Marquardt, we have to solve the system of 
equations~(\ref{260201}) by iterative methods.
It turns out that standard iterative solvers need many iterations 
since the systems becomes very ill-conditioned as
$\gamma_{\it}\to 0$. 

For the efficient solution of
these linear systems we apply the CG-method and exploit its close connection to
Lanczos' method. The latter method is used to approximately compute eigenpairs 
of~$\bl{G}_{\it}^{\top}\bl{G}_{\it}$ to construct a spectral preconditioner for
the CG-method. Since the eigenvalues~$\lambda_1\ge
\ldots\ge\lambda_\dimX$ of~$\bl{A}_{\it}^{\top}\bl{A}_{\it}$ decay at an
exponential rate, it turns out that the approximations determined by Lanczos'
method are well suited to construct an efficient spectral preconditioner.
Spectral preconditioning is reviewed in \S \ref{sec:spectralprec}.
In  \S \ref{sec:itLanczos} we describe how 
the original method proposed in \cite{hohage:01} can be improved by 
the construction of updates of the preconditioner during the Newton 
iteration.
For a convergence analysis of the IRGNM in combination with the 
discrepancy principle \eqref{discrepancy1} discussed below
we refer to~\cite{langer:07, LaHo:07}. It should be mentioned that
all known convergence
results need some condition restricting the degree of nonlinearity of~$\bl{F}$,
and unfortunately none of these conditions has been verified for the 
electromagnetic medium scattering problem. 

\medskip
An essential element of any iterative regularization method for 
an ill-posed problem is a data-driven choice of the stopping index. 
The most common rule
is Morozov's discrepancy principle~\cite{Mor:66}, which consists in
stopping the Newton iteration at the first index~$\It$ satisfying
\begin{equation}\label{discrepancy1} 
\|\bl{F}(\bl{x}_{\It})-\bl{y}^{\rm obs}\| \le \tau\delta <
\|\bl{F}(\bl{x}_{\it})-\bl{y}^{\rm obs}\|,\qquad 0\le \it<{\It},\quad\tau > 1.
\end{equation} 
The discrepancy principle is also frequently used for random noise setting 
$\delta=\sqrt{\mathbb{E}\|\bl{\epsilon}\|^2}$. However, it is easy to see
that this cannot give good results in the limit $N\to\infty$ 
(see e.g.~\cite{BHM:09}), and this is confirmed in our numerical experiments. 
In \S \ref{sec:Lepskij} we show how a Lepski{\u\i}-type stopping rule
can be implemented efficiently in combination with the regularization
method studied in \S \ref{sec:itLanczos}. 

Finally, in \S \ref{sec:num} we report on some numerical experiments
to demonstrate the efficiency of the methods proposed in this paper.

\section{electromagnetic medium scattering problem}\label{sec:forward}
The propagation of time-harmonic electromagnetic waves in an inhomogeneous, non-magnetic,
isotropic  medium without free charges is described by the time-harmonic Maxwell equations
\begin{subequations}\label{eq:forward}
\begin{equation}\label{eq:Maxwell}
\curl\curl \El(\bfr) -\kappa^2(1-a(\bfr))\El(\bfr) = 0,\qquad \bfr\in\Rset^3
\end{equation}
(see \cite{CK:97}). Here $\El:\Rset^3\to\Cset^3$ describes the space-dependent
part of a time-harmonic electromagnetic field of the form 
$\Re \left (\El(\bfr)e^{-i\omega t}\right )$ with angular frequency $\omega>0$. 
Moreover, $\kappa:= \sqrt{\varepsilon_0\mu_0}\omega$
denotes the wave number, $\varepsilon_0$ the electric permittivity of vacuum,
and $\mu_0$ the magnetic permeability of vacuum. The refractive index of the
medium given by 
\[n(\bfr) = \sqrt{1-a(\bfr)}
\] 
is assumed to be $C^{1,\alpha}$-smooth, real and positive in this paper. 
Moreover, we assume that $\supp a\subset B_{1}=\{\bfr\in\Rset^3:|\bfr|<1\}$. 
Now, given a plane incident wave
\[
\Eli(\bfr) = \Eli(\bfr;\bfd,\bfp)=\exp(-i\kappa \bfr\cdot\bfd)\bfp
\]
with direction $\bfd\in S^2$ and polarization $\bfp\in \Cset^3$ 
such that $\bfp\cdot
\bfd=0$, the forward scattering problem consists in finding a total field 
$\El:\Rset^3\to \Cset^3$ satisfying \eqref{eq:Maxwell} such that the
scattered  field $\Els:=\El-\Eli$ satisfies the Silver-M\"uller radiation condition
\begin{equation}\label{eq:SilverMueller}
\lim_{|\bfr|\to\infty}\left (\curl \Els(\bfr)\times \bfr-i\kappa|\bfr|\Els(\bfr)\right )=0
\end{equation}
\end{subequations}
uniformly for all directions $\hat{\bfr}=\bfr/|\bfr|\in S^2$. The latter condition
implies that $\Els$ has the asymptotic behavior 
\[
\Els(\bfr;\bfd,\bfp) = \frac{e^{i\kappa |\bfr|}}{|\bfr|}\Elinf(\hat{\bfr};\bfd,\bfp) + O\left
(\frac{1}{|\bfr|^2}\right ),\qquad |\bfr|\to\infty
\]
with a function $\Elinf(\cdot;\bfd,\bfp):S^2\to\Cset^3$ called the \emph{far field pattern} 
of $\Els$. It satisfies $\Elinf(\hat{\bfr};\bfd,\bfp)\cdot\hat{\bfr}=0$.  

The inverse problem studied in this paper is to reconstruct $a$ given measurements of
$\Elinf(\hat{\bfr};\bfd,\bfp)$ for all $\hat{\bfr},\bfd\in S^2$ and
$\bfp\in \Cset^3$ such that $\bfd\cdot \bfp=0$.

The forward scattering problem has an equivalent formulation in terms of the 
electromagnetic Lippmann-Schwinger equation 
\begin{equation}\label{eq:LippmannSchwinger}
\begin{aligned}
&\El(\bfr )  + \kappa^2\int_{B_1}\!\! \Phi(\bfr -\tilde{\bfr})a(\tilde{\bfr}) 
\El(\tilde{\bfr})\D \tilde{\bfr} 
+ \grad \!\!\int_{B_1}\!\! \Phi(\bfr -\tilde{\bfr}) \frac{\grad a(\tilde{\bfr})}{1-a(\tilde{\bfr})}\cdot \El(\tilde{\bfr}) \D \tilde{\bfr} 
= \El^i(\bfr )
\end{aligned}
\end{equation}
for $\bfr\in\Rset^3$ with the scalar fundamental solution
$\Phi(\bfr):= \exp(i\kappa|\bfr|)/(4\pi |\bfr|)$. 
For the numerical solution of the forward scattering problems 
we use a fast solver of \eqref{eq:LippmannSchwinger}, which converges super-linearly
for smooth refractive indices (see \cite{hohage:06}).

We typically use between $3\cdot 32^3 = 98304$ and
$3\cdot 64^3= 786\,432$ degrees of freedom to represent 
$\El(\cdot;\bfd,\bfp)$ for each $\bfd,\bfp\in S^2$. The unknown perturbation $a$ of the
refractive index is represented by a set of coefficients $\bfx\in\Rset^\dimX$ with 
$500\leq \dimX\leq 2\,000$ using tensor
products of splines in radial direction and spherical harmonics in angular direction
(see \cite{hohage:01}). Moreover, we use 25 incident waves
with random incident directions $\bfd_j$ 
and random polarizations $\bfp_j$ where the directions $\bfd_j$ were drawn
from the uniform distribution on $S^2$.
The exact data are given by complex numbers 
$\Elinf(-\tilde{\bfd}_l;\bfd_j,\bfp_j)\cdot\tilde{\bfp}_l$ 
for $j\in\{1,\dots,25\}$ and 
$l\in\{1,\dots,100\}$ where the $\tilde{\bfd}_j$ and $\tilde{\bfp}_j$ were
generated in the same way as the $\bfd_j$ and $\bfp_j$. This
yields a real data vector $\bfy\in\Rset^{\dimY}$ of size 
$\dimY=2\cdot 25\cdot 100 = 5000$. 

\section{spectral preconditioning}\label{sec:spectralprec}
\subsection{CG method and Lanczos' method}
Let us start by recalling the preconditioned conjugate gradient (CG) method and 
its connection to Lanczos' method (see e.g.~\cite{demmel:97,  GvL:83,
vdV:03}). We consider a preconditioned equation
\begin{equation}\label{040109}
\M^{-1}\bl{G}^{\top}\bl{G}\bl{h} = \M^{-1}\bl{G}^{\top}\bl{g},
\end{equation}
where~$\bl{G}\in \Rset^{\dimY\times\dimX}$ is an arbitrary matrix of rank $\dimX$, and
$\M\in \Rset^{\dimX\times \dimX}$ is a symmetric and positive definite preconditioning
matrix. Although the matrix $\M^{-1}\bl{G}^{\top}\bl{G}$ is not
symmetric in general, the induced linear mapping in $\Rset^{\dimX}$ 
is symmetric and positive definite with respect to the scalar
product~$\left <x,y\right >_{\M}:= \left <\M x,y\right >$ since
\begin{eqnarray*}
\left <\M^{-1}\G^{\top}\G\bl{x},\bl{y}\right >_{\M} & = & \left
<\bl{x},\G^{\top}\G\bl{y}\right > =
\left <\bl{x},\M^{-1}\G^{\top}\G\bl{y}\right >_{\M},\\
\left <\M^{-1}\G^{\top}\G\bl{x},\bl{x}\right >_{\M} & = &
\left < \bl{x},\G^{\top}\G\bl{x}\right > = \|\G\bl{x}\|^2 > 0\quad \mbox{for }\bl{x}\neq 0.
\end{eqnarray*}
Therefore, the CG-method applied to~(\ref{040109}) can be coded as follows:
\begin{alg}\label{211107}\rm (Preconditioned conjugate gradient method)
\small
\begin{itemize}
\item[] $\bl{h}^0 = 0;\hspace{0.1cm} \bl{d}^0 = \bl{g};\hspace{0.1cm}
  \bl{r}^0=\bl{G}^{\top}\bl{d}^0;\hspace{0.1cm} \bl{p}^1=\bl{z}^0=\M^{-1}\bl{r}^0;\hspace{0.1cm} \cgit = 0;$ 
\item[] while $\|\bl{r}^\cgit\| > \varepsilon\gamma\|\bl{h}^{\cgit}\|$
\begin{itemize}
\item[] $\cgit=\cgit+1;$
\item[] $\bl{q}^\cgit = \bl{G}\bl{p}^\cgit;$  
\item[] $\alpha_\cgit = \left <\bl{r}^{\cgit-1}, \bl{z}^{\cgit-1}\right >/\|\bl{q}^\cgit\|^2;$
\item[] $\bl{h}^\cgit = \bl{h}^{\cgit-1}+\alpha_\cgit\bl{p}^\cgit;$ 
\item[] $\bl{d}^\cgit=\bl{d}^{\cgit-1}-\alpha_\cgit\bl{q}^\cgit;$
\item[] $\bl{r}^\cgit=\bl{G}^{\top}\bl{d}^\cgit;$
\item[] $\bl{z}^\cgit=\M^{-1}\bl{r}^{\cgit};$
\item[] $\beta_\cgit = \left <\bl{r}^\cgit,\bl{z}^\cgit\right >/\left <\bl{r}^{\cgit-1}, \bl{z}^{\cgit-1}\right >;$
\item[] $\bl{p}^{\cgit+1} = \bl{z}^\cgit+\beta_\cgit\bl{p}^\cgit$. 
\end{itemize}
\end{itemize}\normalsize
\end{alg}
The stopping criterion~$\|\bl{r}^\cgit\| > \varepsilon\gamma\|\bl{h}^{\cgit}\|$
ensures a relative accuracy of $\varepsilon/(1-\varepsilon)$ of the approximate solution
if $\|((\bf{G}^{\top}\bf{G})^{-1}\| \leq 1/\gamma$ (i.e.~$\gamma=\gamma_{\it}$ if
$\bl{G}=\bl{G}_{\it}$, see \cite{langer:07, LaHo:0801}). 

Quantities arising in Algorithm~\ref{211107} can
be used to approximate the largest eigenvalues and corresponding eigenvectors
of~$\M^{-1}\bl{G}^{\top}\bl{G}$ as follows: Multiplying~$\bl{z}^j =
\bl{p}^{j+1}-\beta\bl{p}^j$ from the left
by~$\|\bl{z}^j\|^{-1}_{\M}\bl{G}$ and using the definitions and identities
\[
\tilde{\bl{z}}^j:=\frac{\bl{z}^j}{\sqrt{(\bl{z}^j)^{\top}\M\bl{z}^j}}, \quad
\tilde{\bl{q}}^j:=\frac{\bl{q}^j}{\|\bl{q}^j\|},\qquad
\frac{\|\bl{q}^{j+1}\|}{\|\bl{z}^j\|_{\M}}=\frac{1}{\sqrt{\alpha_{j+1}}},\quad
\frac{\|\bl{q}^j\|}{\|\bl{z}^j\|_{\M}} = \frac{1}{\sqrt{\alpha_j\beta_j}}
\] 
yields 
\begin{equation}\label{2009-08-25_01}
\bl{G}\tilde{\bl{z}}^0 = \frac{1}{\sqrt{\alpha_1}}\tilde{\bl{q}}^{j},\quad
\bl{G}\tilde{\bl{z}}^j = \frac{1}{\sqrt{\alpha_{j+1}}}\tilde{\bl{q}}^{j+1} -
\sqrt{\frac{\beta_j}{\alpha_j}}\tilde{\bl{q}}^j,\qquad j = 1,\ldots,\cgit-1.
\end{equation}
The identity~$\alpha_j\bl{q}^j = \bl{d}^{j-1}-\bl{d}^j$ multiplied from the left
by~$(\|\bl{q}^j\|\alpha_j)^{-1}\bl{G}^{\top}$ together with 
\[
\frac{\|\bl{z}^{j-1}\|_{\M}}{\|\bl{q}^j\|\alpha_j} =
\frac{1}{\sqrt{\alpha_j}}, \qquad
\frac{\|\bl{z}^{j}\|_{\M}}{\|\bl{q}^j\|\alpha_j} = \sqrt{\frac{\beta_j}{\alpha_j}}
\]
yields
\begin{equation}\label{2009-08-25_02}
\M^{-1}\bl{G}^{\top}\bl{\tilde{q}}^j = \frac{1}{\sqrt{\alpha_j}}\bl{z}^{j-1}
- \sqrt{\frac{\beta_j}{\alpha_j}}\bl{z}^{j},\qquad j = 1,\ldots,\cgit.
\end{equation}
Putting~(\ref{2009-08-25_01}) and~(\ref{2009-08-25_02}) together we have for
all~$j = 1,\ldots,\cgit-1$
\begin{eqnarray*}
\M^{-1}\bl{G}^{\top}\bl{G}\tilde{\bl{z}}^0 & = &
\frac{1}{\alpha_1}\tilde{\bl{z}}^0 -
\frac{\sqrt{\beta_1}}{\alpha_1}\tilde{\bl{z}}^1,\\
\M^{-1}\bl{G}^{\top}\bl{G}\tilde{\bl{z}}^j & = & 
-\frac{\sqrt{\beta_j}}{\alpha_j}\tilde{\bl{z}}^{j-1} + \left
 (\frac{1}{\alpha_{j+1}} + \frac{\beta_j}{\alpha_j}\right )\tilde{\bl{z}}^{j} - 
\frac{\sqrt{\beta_{j+1}}}{\alpha_{j+1}} \tilde{\bl{z}}^{j+1}.
\end{eqnarray*}
These formulas can be rewritten as
\begin{equation}\label{2009-08-25_03}
\M^{-1}\bl{G}^{\top}\bl{G}\bl{Z}_\cgit = \bl{Z}_\cgit\bl{T}_\cgit-\left
(0,\ldots,0,\frac{\sqrt{\beta_\cgit}}{\alpha_\cgit}\tilde{\bl{z}}^\cgit\right )
\end{equation}
where~$\bl{Z}_\cgit:=(\tilde{\bl{z}}^0,\ldots,\tilde{\bl{z}}^{\cgit-1})$ and
\[
\bl{T}_\cgit := \left (
\begin{array}{ccccc}
\frac{1}{\alpha_1} & -\frac{\sqrt{\beta_1}}{\alpha_1} \\
-\frac{\sqrt{\beta_1}}{\alpha_1} & \frac{1}{\alpha_2}+\frac{\beta_1}{\alpha_1} & 
-\frac{\sqrt{\beta_2}}{\alpha_2} \\
& \ddots & \ddots & \ddots \\
 & & \ddots &\ddots & -\frac{\sqrt{\beta_{\cgit-1}}}{\alpha_{\cgit-1}}\\
& & & -\frac{\sqrt{\beta_{\cgit-1}}}{\alpha_{\cgit-1}} & \frac{1}{\alpha_{\cgit}}+\frac{\beta_{\cgit-1}}{\alpha_{\cgit-1}}.
\end{array}\right ),
\]
If we denote by~$\theta_1>\ldots>\theta_\cgit>0$ and~$\bl{v}_1,\ldots,\bl{v}_\cgit$ 
the eigenvalues with corresponding eigenvectors of the symmetric and positive
definite matrix~$\bl{T}_\cgit$, \eqref{2009-08-25_03} implies
\[
\bl{Z}_\cgit^{\top}\M^{-1}\G^{\top}\G\bl{Z}_{\cgit}\bl{v}_j = \theta_j\bl{v}_j,\qquad j
= 1,\ldots,\cgit.
\]
Hence, in the case that~$\bl{z}^\cgit$ vanishes~$\theta_1>\ldots>\theta_\cgit$ are exact
eigenvalues of~$\M^{-1}\bl{G}^{\top}\bl{G}$
with corresponding eigenvectors~$\bl{Z}_\cgit\bl{v}_1,\ldots,\bl{Z}_\cgit\bl{v}_\cgit$. 
In the typical case~$\bl{z}^\cgit\neq 0$ the
vectors~$\bl{Z}_\cgit\bl{v}_j$ usually converge rapidly to the eigenvectors
corresponding to the outliers in the spectrum
of~$\M^{-1}\bl{G}^{\top}\bl{G}$ (cf.~\cite{demmel:97, GvL:83} and the
references on the Kaniel-Paige theory therein) and Lanczos' method can be
interpreted as a particular case of the Rayleigh-Ritz method. This connection
can be used to interpret the so-called \emph{Ritz values}~$\theta_1>\ldots
>\theta_\cgit$ and the \emph{Ritz vectors}~$\bl{Z}_\cgit\bl{v}_1,\ldots,\bl{Z}_\cgit\bl{v}_\cgit$ 
as approximations to some eigenpairs of~$\M^{-1}\G^{\top}\G$. 

If $\bl{W}\bl{\Lambda}\bl{W}^{\top}$ is an eigendecomposition of the matrix
$\bl{T}_\cgit$ with $\bl{M}=\bl{I}$, i.e.
$\bl{W}=(\bl{w}_1,\ldots,\bl{w}_{\tilde{\cgit}})$ is orthogonal and
$\bl{\Lambda}=\mbox{diag}(\theta_1^{(2)},\ldots,\theta_{\tilde{\cgit}}^{(2)})$,
one can prove the equality~(see~\cite{langer:07})
\begin{equation}\label{eq:accuracyRitz}
\|\G^{\top}\G(\bl{Z}_{\tilde{\cgit}}\bl{w}_i)-(\bl{Z}_{\tilde{\cgit}}\bl{w}_i)\theta_i^{(2)}\| =
\frac{\sqrt{\beta_{\tilde{\cgit}}}}{\alpha_{\tilde{\cgit}}}|\bl{w}_i(\tilde{\cgit})|,\qquad
i = 1,\ldots,\tilde{\cgit},
\end{equation}
where~$\bl{w}_i(\tilde{\cgit})$ denotes the bottom entry of~$\bl{w}_i$. 
This identity can be used to judge the accuracy of the Ritz pairs and to decide
which of them to use in the spectral preconditioner.

\subsection{Spectral preconditioning with Tikhonov regularization}\label{090301}
Assume now that $\G$ is of the special form
\[
\G := \left (\begin{array}{c} 
\bl{A}\\ 
\sqrt{\gamma}\;\bl{I}\end{array}\right ) 
\]
with $\bl{A}\in\Rset^{\dimY\times \dimX}$. 
Let $\bl{u}_1,\dots,\bl{u}_\dimX$ be orthonormal
eigenvectors of $\bl{A}^{\top}\bl{A}$, and let $\lambda_1,\dots,
\lambda_\dimX$ be the corresponding eigenvalues. 

Given eigenpairs $(\lambda_j,\bl{u}_j)$ for $j$ in some non-empty subset
$\calJ\subset \{1,2,\dots,\dimX\}$,  we define a spectral preconditioner
for $\bl{G}^{\top}\bl{G}=\gamma \bl{I}+ \bl{A}^{\top}\bl{A}$ by
\[
\M := \gamma\bl{I} + \sum_{j\in \calJ}\lambda_j
\bl{u}_j(\bl{u}_j)^{\top}.
\]
Its properties are summarized in the following proposition:

\begin{prop}\label{prop:precond}
Assume that $\mathrm{rank}(\mathbf{A})=\dimX$. Then
\begin{itemize}
\item[a)] $\M$ is symmetric and positive definite, 
and its inverse is given by 
\[
\M^{-1} = \frac{1}{\gamma}\bl{I}+\sum_{j\in\calJ}
\left(\frac{1}{\lambda_j+\gamma}-\frac{1}{\gamma}\right )\bl{u}_j\bl{u}_j^{\top}. 
\] 
\item[b)] $\M^{-1}\G^{\top}\G \bl{x}
= \bl{x} +
\sum_{j\notin\calJ} \frac{\lambda_j}{\gamma}
\left <\bl{x},\bl{u}_j\right >\bl{u}_j
=\G^{\top}\G\M^{-1}\bl{x}$. 
\item[c)] The spectrum of the preconditioned matrix is given by
\begin{equation}\label{specass}
\sigma(\M^{-1}\G^{\top}\G) =
 \left\{1+\frac{\lambda_{j}}{\gamma}:j\notin\calJ\right\} \cup \{1\},
\end{equation}
and the eigenvalue $1$ has multiplicity $\# \calJ$.
\item[d)] If $\mu\neq 1$ is an eigenvalue of $\M^{-1}\G^T\G$
with corresponding eigenvector $\mathbf{u}$,
then $(\gamma(\mu-1),\mathbf{u})$ is an eigenpair of $\mathbf{A}^{\top}\mathbf{A}$.
\end{itemize} 
\end{prop}

\pro
$\M$ is obviously symmetric, and it is positive definite
since all its eigenvalues are $\geq\gamma>0$. The formula for the inverse follows from a straightforward computation. 

Let $\mathcal{U}:=\mathrm{span}\{u_j:j\in \calJ\}$. Identifying matrices with their induced
linear mappings, we have $\M|_{\mathcal{U}} = \G^{\top}\G|_{\mathcal{U}}$ 
and $\M|_{\mathcal{U}^{\perp}} = \gamma\bl{I}|_{\mathcal{U}^{\perp}}$,
and $\mathcal{U}$ and $\mathcal{U}^{\perp}$ are invariant under all the involved linear mappings. 
Therefore, 
\begin{eqnarray*}
&& \M^{-1}\G^{\top}\G|_{\mathcal{U}} = \bl{I}|_{\mathcal{U}}
= \G^{\top}\G \M^{-1} |_{\mathcal{U}},\\
&& \M^{-1}\G^{\top}\G|_{\mathcal{U}^{\perp}} 
= \frac{1}{\gamma} \G^{\top}\G|_{\mathcal{U}^{\perp}}
= \G^{\top}\G\M^{-1} |_{\mathcal{U}^{\perp}}.
\end{eqnarray*}
Since $\G^{\top}\G\bl{x} = \sum_{j=1}^M(\gamma+\lambda_j)\left<\bl{x},\bl{u}_j\right> \bl{u}_j$ 
for all $\bl{x}\in\Rset^M$, assertion b) follows. 
c) is obtained from b) by inserting the eigenvectors~$\bl{u}_j$ into the formula.

If $(\mu,{\bf u})$ is an eigenpair of $\M^{-1}\G^T\G$
and $\mu\neq 1$, it follows that $\bl{u}\in \mathcal{U}^{\perp}$. Therefore,
$\mu \bl{u} = \frac{1}{\gamma} \G^T\G\bl{u} 
= \bl{u} + \frac{1}{\gamma} \bl{A}^{\top}\bl{A}\bl{u}$.
This implies assertion d).
%
\endpro

\begin{rem}
We comment on the assumption $\mathrm{rank}(\mathbf{A})=\dimX$ in
Proposition \ref{prop:precond}. For the acoustic medium
scattering problem injectivity of the continuous Fr\'echet derivative
$F'[\bl{x}]$ has been shown in \cite[Prop.~2.2]{hohage:01}. 
For the electromagnetic medium scattering problem uniqueness proofs
for the nonlinear inverse problem (see \cite{CP:92,haehner:00}) 
can be modified analogously
to show injectivity of $F'[\bl{x}]$. It is easy to see that this implies
injectivity of $\bl{A}$ if $\bl{A}$ is a sufficiently accurate discrete
approximation of $F'[\bl{x}]$ on a finite dimensional subspace.
\end{rem}

\section{IRGNM with updated spectral preconditioners}\label{sec:itLanczos}


Spectral preconditioning in Newton methods 
is particularly useful for exponentially ill-posed problems
such as the electromagnetic inverse medium scattering problem.
Typically,  Lanczos' method approximates
  outliers in the spectrum well, whereas eigenvalues in the
  bulk of the spectrum are harder to approximate. Frequently
  the more isolated an eigenvalue is, the better the approximation 
  (see~\cite{kuijlaars:00} and~\cite[Chapter 7]{demmel:97}). For exponentially
  ill-conditioned problems the spectrum of~$\G_{\it,m}^{\top}\G_{\it,m}$ consists of a small
  number of isolated eigenvalues and a large number of eigenvalues clustering at~$\gamma_{\it}$. 
  If all the large isolated eigenvalues are found and computed accurately, spectral
  preconditioning reduces the condition number significantly. 

Updating the preconditioner may be necessary for the following reasons:
\begin{itemize}
\item If the matrix $\G_{\it,m}^{\top}\G_{\it,m}$ has multiple isolated eigenvalues, the Lanczos'
  method approximates at most one Ritz pair corresponding to this multiple eigenvalue. 
\item During Newton's method the regularization parameter~$\gamma_{\it}$ tends
  to zero. Hence, if we keep the number of known eigenpairs for the construction of the
  preconditioner~$\M_{\it,m}$ fixed, the number of CG-steps will increase rapidly
  during our frozen Newton method (see~\cite{LaHo:0801}).
\end{itemize}


In the preconditioned Newton iteration we keep the Jacobian $\bl{A}_m={\bf F}'(\bl{x}_m)$ frozen for
several Newton steps and replace eq.~\eqref{040109} by
\begin{subequations}
\begin{equation}\label{18010802}
\M_{\it,m}^{-1}\G_{\it,m}^{\top}\G_{\it,m}\bl{h} =
\M_{\it,m}^{-1}\G_{\it,m}^{\top}\bl{g}_{\it},
\end{equation}
where
\begin{equation}\label{eq:defiGg}
\G_{\it,m} := \left (\begin{array}{c} 
\bl{A}_{m}\\ 
\sqrt{\gamma_{\it}}\;\bl{I}\end{array}\right ) 
\qquad\mbox{and}\qquad
\bl{g}_{\it} := \left (\begin{array}{c}
\bl{y}^{\rm obs} - \bl{F}(\bl{x}_{\it})\\ 
\sqrt{\gamma_{\it}}\bl{b}_{\it} \end{array}
\right). 
\end{equation}
\end{subequations}
Moreover, given some eigenpairs $\{(\lambda_j^{(m)},\bl{u}_j^{(m)}):j\in\calJ\}$ 
of $\bl{A}_m^{\top}\bl{A}_m$ with orthonormal eigenvectors $\bl{u}_j^{(m)}$,
the spectral preconditioner is defined by
\begin{equation}\label{eq:spectral_prec}
\M_{\it,m} := \gamma_{\it}\bl{I} + \sum_{j\in \calJ}\lambda_j^{(m)}
\bl{u}_j^{(m)}(\bl{u}_j^{(m)})^{\top}.
\end{equation}

A preconditioned semi-frozen Newton method with updates of the preconditioner 
can be coded as follows:
\begin{alg}\label{alg:IRGNM_Prec_Update} {\rm Input: initial guess $\bl{x}_0$,
data $\bl{y}^{\delta}$, $\delta$ and/or $\mathbf{Cov}_{\epsilon}$
\textit{(see \eqref{eq:noise_model})}}
\begin{itemize}\rm
\item[] $\it:=0$; $m:=0$
\item[] repeat
   \begin{itemize}
    \item[+] Evaluate $\bl{F}(\bl{x}_\it)$ and define
     $\bl{G}_{\it,m}, \bl{g}_\it$ by \eqref{eq:defiGg}; 
   \item[] if $\sqrt{\it+1}\geq \sqrt{m+1}+1$
      \begin{itemize}
      \item[$\bullet$] $m:=\it$;
      \item[$\bullet$] Solve $\G_{\it,\it}^{\top}\G_{\it,\it}\bl{h}_\it = \G_{\it,\it}^{\top}\bl{g}_\it$ by CG-method;
      \item[$\bullet$] Compute via Lanczos' method orthonormal Ritz pairs 
        $\{(\mu_j^{(m)},\bl{u}_j^{(m)}):j\in\tilde{\calJ}\}$ of $\G_{k,k}^{\top}\G_{k,k}$;
      \item[$\bullet$] Select subset $\calJ\subset \tilde{\calJ}$
        \textit{(see \eqref{eq:accuracyRitz})}
        and set $\lambda_j^{(m)}:=\gamma_{\it}(\mu_j^{(m)}-1)$ for $j\in\calJ$
         \textit{(see Prop. \ref{prop:precond})};
      \end{itemize}
   \item[] else
      \begin{itemize}
      \item[$\bullet$] Define $\M_{\it,m}$ by \eqref{eq:spectral_prec};
      \item[] if MustUpdate() \hspace*{1cm} \textit{(see Remark \ref{it:MustUpdate} below!)}
       \begin{itemize}
      \item[$\to$] Solve 
              $\M_{\it,m}^{-1/2}\G_{\it,m}^{\top}\G_{\it,m}\M_{\it,m}^{-1/2} 
              \bl{\tilde{h}}_{\it} =
      \M_{\it,m}^{-1/2}\G_{\it,m}^{\top}\bl{g}_{\it}$ by CG-method; 
      \item[$\to$] $\bl{h}_\it:= \M_{\it,m}^{1/2}\bl{\tilde{h}}_{\it}$;
      \item[$\to$] Compute Ritz pairs 
        $\{(\mu_j^{(m)},\bl{u}_j^{(m)}):j\in\tilde{\calJ}_2\}$
        of $\M_{\it,m}^{-1/2}\G_{\it,m}^{\top}\G_{\it,m}\M_{\it,m}^{-1/2}$
                using Lanczos' method;
      \item[$\to$]  Select subset $\calJ_2\subset \tilde{\calJ}_2$
        \textit{(see Remark \ref{it:Lanczos})} 
        and set $\lambda_j^{(m)}:=\gamma_{\it}(\mu_j^{(m)}-1)$ for $j\in\calJ_2$;
      \item[$\to$] Set $\calJ:=\calJ\cup \calJ_2$ and reorthogonalize
        $\{\bl{u}_j^{(m)}:j\in\calJ\}$;
       \end{itemize}
      \item[] else
      \begin{itemize}
      \item[$\to$] Solve $\M_{\it,m}^{-1}\G_{\it,m}^{\top}\G_{\it,m}\bl{h}_{\it} =
      \M_{\it,m}^{-1}\G_{\it,m}^{\top}\bl{g}_{\it}$ by CG-method;
       \end{itemize}
       \item[] end
      \end{itemize}
   \item[] end
   \item[$+$] $\bl{x}_{\it+1}:=\bl{x}_{\it} + \bl{h}_{\it}$; $\it:=\it+1$;
   \end{itemize}
\item[] until Stop() \hspace{1ex} \textit{(see \S \ref{sec:Lepskij})}
\item[] Select stopping index $\It$ 
\textit{(see \S \ref{sec:Lepskij})}
and return $\bl{x}_{\It}$; 
\end{itemize}
\end{alg}

We add some remarks on heuristics and implementation details for 
Algorithm \ref{alg:IRGNM_Prec_Update}:
\begin{enumerate}
\item Usually round-off errors cause loss of orthogonality in the residual 
  vectors~$\bl{z}^j$ computed in Algorithm~\ref{211107}. This
  loss of orthogonality is closely related to the convergence of the Ritz
  vectors (see~\cite{demmel:97,langer:07}). To sustain stability,
  Algorithm~\ref{211107} was amended by 
  a complete reorthogonalization scheme based on Householder transformations 
  (see \cite{GvL:83}).
\item The necessity of reorthogonalization is also our reason for preconditioning
with $\bl{M}_{\it,m}^{-1/2}$ from both sides instead of $\bl{M}_{\it,m}^{-1}$ from the left
when updating the preconditioner. In the latter case, reorthogonalization would have
to be performed with respect to the inner product $\langle\cdot,\cdot\rangle_{\bl{M}_{\it,m}}$,
which is more complicated. Note that
\begin{eqnarray*}
\bl{M}_{\it,m}^{-1/2}\bl{x} &=& \frac{1}{\sqrt{\gamma_{\it}}} \bl{x}
+ \sum_{j\in\calJ}\left(\frac{1}{\sqrt{\gamma_\it + \lambda_j}}
-\frac{1}{\sqrt{\gamma_\it}}\right) (\bl{u}_j^{\top}\bl{x}) \bl{u}_j,\\
\bl{M}_{\it,m}^{1/2}\bl{x} &=& \sqrt{\gamma_{\it}} \bl{x}
+ \sum_{j\in\calJ}\left(\sqrt{\gamma_\it + \lambda_j}
-\sqrt{\gamma_\it}\right) (\bl{u}_j^{\top}\bl{x}) \bl{u}_j.
\end{eqnarray*}
\item Spectrally preconditioned linear systems react very sensitively to errors
  in the eigenelements (see~\cite{GiGr:06, langer:07}). Hence, to ensure
  efficiency of the preconditioner it is necessary that the approximations of
  the Ritz pairs used in the construction of the preconditioner be of high 
  accuracy. This is achieved by choosing $\varepsilon =10^{-9}$ in Algorithm~\ref{211107}
  when updating or recomputing the preconditioner, whereas $\varepsilon=1/3$ is
  sufficient otherwise. Numerical experience 
  shows that computation time invested into improved accuracy
  of the Ritz pairs pays off in the following Newton steps.
\item \label{it:MustUpdate} MustUpdate(): We update the preconditioner if the
last update or recomputation is at least 4 Newton steps ago and the number of inner
iterations in the previous Newton step is $>5$. 
\item We found it useful not to perform a complete recomputation of the current preconditioner 
if it works well. Therefore, we amend the condition $\sqrt{\it+1}\geq \sqrt{m+1}+1$ by the additional requirement that the number of inner iterations in the previous step be not too small, say $>8$.
The condition  $\sqrt{\it+1}\geq \sqrt{m+1}+1$ is a generalization of the rule
to recompute the preconditioner whenever $\it+1$ is a square number, which
was proposed in the original paper \cite{hohage:01}. Under certain conditions it 
was shown in \cite{LaHo:0801} to be optimal among all rules where $\sqrt{\cdot}$ is replaced by a 
function $x\mapsto x^{\mu}$ with $\mu\in(0,1]$.  
\item\label{it:Lanczos} For updating the preconditioner we only select Ritz values of 
$\M_{\it,m}^{-1/2}\G_{\it,m}^{\top}\G_{\it,m}\M_{\it,m}^{-1/2}$ which are sufficiently
well separated from the cluster at $1$, say $\geq 1.1$. First, these eigenvalues
are usually computed more accurately by Lanczos' method, and second, they are more
relevant for preconditioning. 
\item In the initial phase when the updates $\bl{h}_\it$ are large, keeping the
Jacobian frozen is not efficient. Therefore, we use other methods in this phase,
e.g.~Newton-CG. In some cases globalization strategies will be necessary in this phase, 
although this was not the case in the examples reported below. 
\end{enumerate}

\section{Implementation of a Lepski{\u\i}-type stopping rule}
\label{sec:Lepskij}
Lepski{\u\i}-type stopping rules for regularized Newton methods have been studied
in \cite{BH:05,BHM:09}. We refer to the original paper \cite{lepskij:90} on
regression problems and to \cite{MP:03a,mathe:06} for a considerable simplification
of the idea and its application to linear inverse problems. 
As opposed to the discrepancy principle, Lepski{\u\i}-type stopping rules yield
order optimal rates of convergence for all smoothness classes up to the
qualification of the underlying linear regularization method (in case of
random noise typically only up to a logarithmic factor). 

A crucial element of Lepski{\u\i}'s method are estimates of the propagated
data noise error, and the performance depends essentially on the sharpness of these
estimates. Let ${\bf R}_{\it}:=({\bf G}_{\it}^{\top}{\bf G}_{\it})^{-1}{\bf A}_{\it}$. 
If $\mathbf{\epsilon}\in \mathbb{R}^\dimY$ is a deterministic noise vector,
an estimate of the propagated data noise error is given by
\begin{equation}\label{eq:det_noise_estim}
\|{\bf R}_{\it}\mathbf{\epsilon}\|\leq \|{\bf R}_{\it}\|\,\|\mathbf{\epsilon}\|  \leq \frac{1}{2\gamma_{\it}} \|\mathbf{\epsilon}\|,
\end{equation}
and these estimates are sharp if $(\gamma_{\it},\mathbf{\epsilon})$ is
an eigenpair of ${\bf A}_{\it}{\bf A}_{\it}^{\top}$. However, if $\mathbf{\epsilon}$ is a random vector
with $\mathbb{E} \mathbf{\epsilon} =0$, finite second moments with covariance
matrix ${\bf Cov}_{\mathbf{\epsilon}} = (\mathrm{Cov}(\epsilon_i,\epsilon_j))_{i,j=1..\dimY}$, 
the estimate \eqref{eq:det_noise_estim} is usually very pessimistic, and we have
\begin{equation}\label{eq:ran_noise_estim}
\sqrt{\mathbb{E} \|{\bf R}_{\it}\mathbf{\epsilon}\|^2} = \sqrt{\mathrm{trace}({\bf R}_{\it}^{\top}{\bf Cov}_{\mathbf{\epsilon}}{\bf R}_{\it})}.
\end{equation}
Denoting the right hand side of \eqref{eq:det_noise_estim} or \eqref{eq:ran_noise_estim},
respectively, by $\Phi(\it)$, the Lepski{\u\i} stopping rule is defined by
\begin{equation}
\label{eq:defi_Lepskij}
\It_{\rm bal} := \min\{\it\leq \It_{\max}: \|{\bf x}_{\it}-{\bf x}_m\|\leq \rho \Phi(m),
m=\it+1,\dots,\It_{\max}\}
\end{equation}
with a parameter $\rho>4$ and a maximal Newton step number $\It_{\max}$.
We choose $\rho=4.1$ in our numerical experiments and
$\It_{\max}:=\max\{\it\in\mathbb{N}:\Phi(\it)\leq R\}$ with a reasonable
upper bound on the size of propagated data noise in the optimal reconstruction.
$R$ may be an a-priori known bound $\|\bf{x}-\bf{x}_0\|$. However, it is advisable
to choose a smaller value of $R$ to reduce the number of Newton iterations. The
final results ${\bf x}_{\It_{\rm bal}}$ do not depend critically on $R$.

The main computational challenge in the implementation of the stopping rule \eqref{eq:defi_Lepskij} 
for random noise is the efficient and accurate computation of  $\Phi(\it)$.  One possibility is
to generate $L$ independent copies $\mathbf{\epsilon}_1,\dots,\mathbf{\epsilon}_L$ of the noise vector and use the approximation $\Phi(\it)\approx (L^{-1}\sum_{l=1}^L \|{\bf R}_{\it}\mathbf{\epsilon}_l\|^2)^{1/2}$. 
However, this involves the iterative solution of $L+1$ instead of $1$ least squares system
and leads to a tremendous increase of the computational cost.

With the methods described in the previous sections we can construct 
approximations ${\bf R}_{\it}^{\rm app}:= \sum_{j\in \calJ_m} 
\frac{\sqrt{\lambda_j}}{\gamma_{\it}+\lambda_j}{\bf u}_j
{\bf w}_j^{\top}$ of ${\bf R}_{\it}$, which allow cheap matrix-vector 
multiplications not involving evaluations of the forward mapping ${\bf F}$. 
This yields the approximation
\begin{equation}
\Phi(k) \approx  \left(\frac{1}{L}\sum_{l=1}^L \|{\bf R}_{\it}^{\rm app}\mathbf{\epsilon}_l\|^2
\right )^{1/2}. 
\label{eq:approx_ran_noise_estim}
\end{equation}
Here ${\bf w}_j:={\bf A}_m {\bf u}_j / \| {\bf A}_m {\bf u}_j\|$ denote the approximated
left singular vectors of ${\bf A}_m$, which can be computed directly by an appropriately modified
Lanczos method (see e.g. \cite{GvL:83}). 
In the case of white noise, i.e. ${\bf Cov}_{\mathbf{\epsilon}} = \sigma^2 {\bf I}_{\dimY}$, the
expected value of the right hand side is given by the simple expression
\begin{equation}
\Phi(\it) \approx \sigma \left(\sum_{j=\calJ_m} \frac{\lambda_j}{(\gamma_{\it}+\lambda_j)^2}\right)^{1/2}
\label{eq:approx_white_noise_estim}.
\end{equation}
Obviously, equality holds in \eqref{eq:approx_white_noise_estim} if 
$(\lambda_j^2)_{j\in \calJ_m}$ is a complete set of eigenvalues of ${\bf A}_m^{\top}{\bf A}_m$
(with multiplicities). Under certain assumptions it has been shown in \cite{BHMR:07}
in an infinite dimensional setting that the left hand side can be bounded
by a small constant times the right hand side uniformly in $\gamma_{\it}$ if
$(\lambda_j^2)_{j\in \calJ_m}$ contains all eigenvalues $\geq \gamma_{\it}$. 
Our numerical results in section \ref{sec:num} indicate that this approximation is sufficiently
accurate. 

\begin{figure}[ht]
\hspace*{-0.2\textwidth}\includegraphics[width=1.3\textwidth]{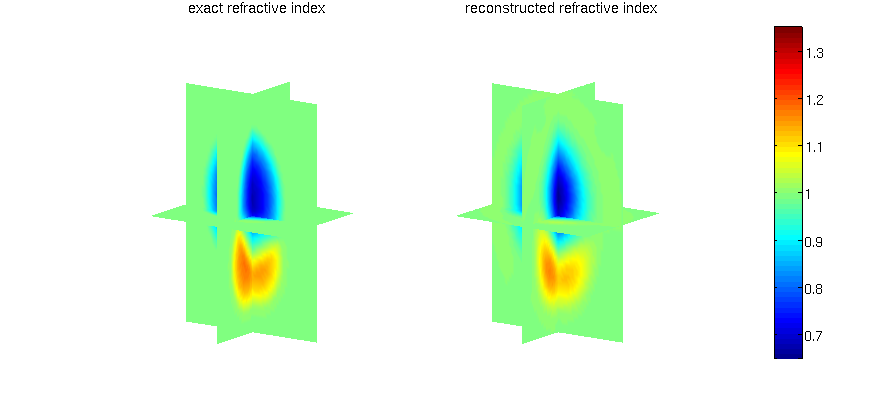}
\caption{\label{fig:two_bumps}
An exact refractive index and its reconstruction by
the IRGNM with updated preconditioner at iteration 23.
The plots show the cube $[-1,1]^3$, the wave number is $\kappa=1$.}
\end{figure}

\section{Numerical results}\label{sec:num}




As a test example we use the refractive index shown in Figure~\ref{fig:two_bumps}.
For further information on the forward problem and its numerical solution we
refer to \S \ref{sec:forward} and \cite{hohage:06}.

\begin{figure}[ht]
\includegraphics[width=\textwidth]{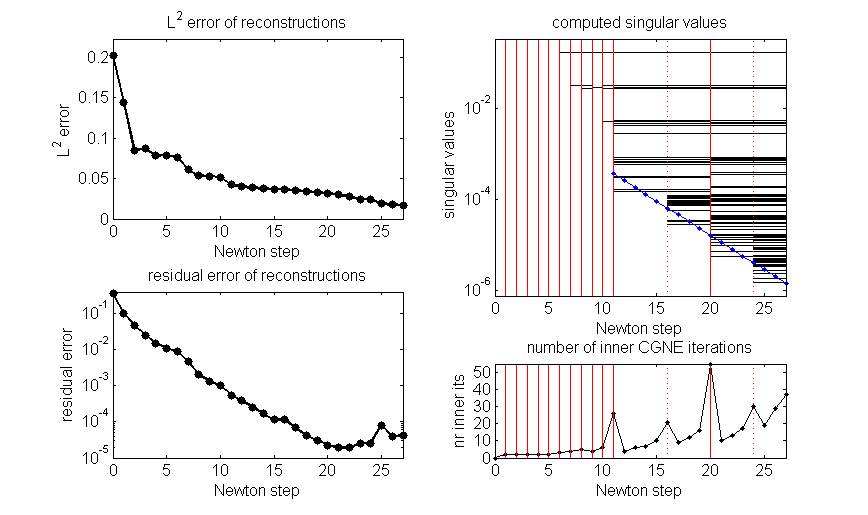}
\caption{\label{fig:two_bumps_Stefan1}Performance of the preconditioned
Newton method applied to the example in Figure \ref{fig:two_bumps}. 
The left panels show the continuous $L^2$-error of the reconstructed
refractive indices and the norm of the residuals $\|{\bf F}({\bf x}_k)-{\bf y}^{\rm obs}\|_2$
over the Newton step $\it$. The upper right panel shows the computed singular values used
for preconditioning (solid horizontal lines). After step 11 the method changed from Newton-CG to IRGNM. 
At steps $m=11$ and $20$ the preconditioner was completely recomputed for the derivative at a new iterate. This is indicated by solid vertical lines. Updates of
the preconditioner, indicated by dotted vertical lines, were performed at steps 16 and 24.
The dots on the diagonal line indicate the values of the regularization parameter
$\sqrt{\gamma_{\it}}$. Some sufficiently accurate singular values computed
during the Newton-CG phase are shown for their own interest although they were not
used in the computation. The right lower panel shows the number of inner CGNE iterations
over the Newton step $\it$}
\end{figure}

Figure \ref{fig:two_bumps_Stefan1} illustrates the performance of the IRGNM with updated preconditioners.
In the update steps for the preconditioner at  $\it=16$ and $24$ the new singular
values mainly fall into two categories: First, we have singular values which are 
not well separated from the cluster for the regularization parameter 
$\gamma_m$ but are well separated for $\gamma_{\it}$. These singular values are in
or near the interval $[\sqrt{\gamma_{\it}},\sqrt{\gamma_m}]$. The second category
are multiple or nearly multiple singular values where only one element in the
eigenspace is found in the application of Lanczos' method. The use of
an update clearly reduces the number of inner CGNE steps in the following Newton 
iterations.

\begin{figure}[ht]
\begin{center}
\includegraphics[width=0.9\textwidth]{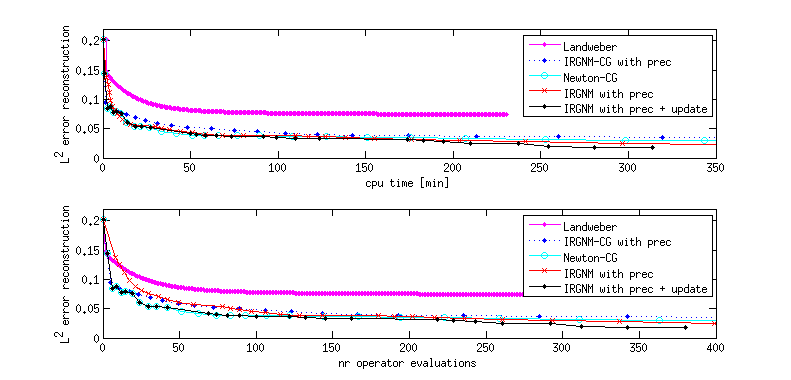}
\end{center}
\caption{\label{fig:comparison}Work-precision diagram comparing the performance
of different methods applied to the problem in Figure \ref{fig:two_bumps} for
exact data.}
\end{figure}

Moreover, in Fig \ref{fig:comparison} 
we compared the speed of convergence of the iterative regularization methods 
discussed in the introduction for exact data. Here we measure  "speed" both in terms
of cpu-time and in terms of the number of evaluations of $\bl{F}$, $\bl{F}'[\bl{x}_m]$ or 
$\bl{F}'[\bl{x}_m]'$. Landweber iteration is clearly the slowest method although some
good progress is achieved in the first few steps. The Newton-CG method
performs very well up to some accuracy after which on it becomes
slow, a behavior also observed in most other examples. We stopped the Newton-CG
iteration at an $L^2$-error of $\approx 0.28$, which was achieved by the updated
preconditioned IRGNM 2.5 times earlier. 
We also include a comparison with the preconditioned IRGNM without updating
as suggested in \cite{hohage:01}. The updating improves the performance particularly
at high accuracies. Note that in
the first Newton steps where $a$ is still small, the iterative solution of 
the forward problem is faster than in later Newton steps. 

\begin{table}
\begin{center}
\begin{tabular}{|l|c|c|}
\hline
stopping rule & stopping index & $L^2$ error at stopping index \\ \hline\hline
optimal  & $15.93 \pm 0.57$ & $0.0406 \pm 0.0010$ \\ \hline
Lepski{\u\i}  & $12.20 \pm 0.54$ & $0.0474 \pm 0.0016$ \\ \hline
discrepancy & $4.87\pm 0.43$ & $0.0744 \pm 0.0001$ \\ \hline
\end{tabular}
\end{center}
\caption{\label{tab:stopRules} Performance of stopping rules averaged over 
15 noise samples for the problem in Figure~\ref{fig:two_bumps}. The numbers
indicate means and standard deviations.}
\end{table}

Finally, we tested the performance of the Lepski{\u\i}-type stopping for
randomly perturbed data. More precisely, we added independent Gaussian variables 
to each data point. The ''relative noise level'' $\|\epsilon\| / \|\bl{y}\|$ was
about $2\%$, but we stress that such a point-wise definition of the noise level does not
make sense for random noise when considering the limit $N\to\infty$. We compare
the discrepancy principle with $\tau=2$ to Lepski{\u\i}'s method with $\rho=4.1$.
Moreover, we look at the optimal stopping index for each noise sample.
As expected, the discrepancy principle stops the iteration too early.
Note in Figure~\ref{fig:two_bumps_Stefan1} that $\|\bl{F}(\bl{x}_\it)-\bl{F}(\bl{x})\|$ is at least
an order of magnitude smaller than $\|\varepsilon\|$ at the optimal $k\approx 16$.
(In Figure~\ref{fig:two_bumps_Stefan1} we used exact data, but the behavior is similar
for noisy data.) 
The results in Tab.~\ref{tab:stopRules} indicate that Lepski{\u\i}'s stopping rule is stable and
yields considerably better results than the discrepancy principle.


\paragraph*{Acknowledgments:}
The second author acknowledges financial support by DFG
(German Research Foundation) in the Research Training Group 1023
''Identification in Mathematical Models: Synergy of Stochastic and Numerical
Methods''. 


\end{document}